\documentclass[12pt,onecolumn,leqno]{article}

\usepackage{graphicx}
\usepackage{indentfirst}
\usepackage{amsmath}
\usepackage{amsfonts}
\usepackage{amssymb}
\usepackage{color}
\usepackage[mathscr]{eucal}
\usepackage{latexsym}

\newtheorem{theorem}{Theorem}[section]
\newtheorem{lemma}{Lemma}[section]
\newtheorem{defini}{Definition}[section]

\newtheorem{prop}{Proposition}[section]
\newtheorem{rem}{Remark}[section]
\numberwithin{equation}{section}

\title{BC$_{{}_{\text N}}$-graded Lie algebras arising from fermionic representations
} \author{Hongjia Chen  and  Yun Gao\footnote{Research was
partially supported by NSERC of Canada and Chinese Academy of
Science.} }
    \date{}  % blank {} to invalidate the current date

\begin{document}

%-------------- customize your commands
%\newcommand{\supercite}[1]{\textsuperscript{\cite{#1}}}

%-------------- the title
\maketitle              %生成标题

     \begin{abstract}
         We use fermionic representations to obtain a class of
          BC$_{{}_{\text N}}$-graded Lie algebras coordinatized by quantum tori with  nontrivial central extensions.

     \end{abstract}

%-------------- main part
%\part{Introduction}    %Part II ……，部分命令

\setcounter{section}{-1}       %给章节赋值，使得Introduction编号为0

\section{Introduction}
%\vspace{0.2cm}

Lie algebras graded by the reduced finite root systems were first
introduced by Berman-Moody [BM] in order to understand the
generalized intersection matrix algebras of Slodowy. [BM]
classified
 Lie algebras graded by the root systems of type $A_l, l\geq 2$,
$D_l, l\geq 4$ and $E_6, E_7, E_8$ up to central extensions.
Benkart-Zelmanov [BZ] classified Lie algebras graded by the root
systems of type $A_1, B_l, l\geq 2$, $C_l, l\geq 3$, $F_4$ and
$G_2$ up to central extensions. Neher [N] gave a different
approach for all reduced root systems except $E_8, F_4$ and $G_2$.
The idea of root graded Lie algebras can be traced back to Tits
[T] and Seligman [S]. [ABG1] completed the classification of the
above root graded Lie algebras by figuring out explicitly the
centers of the universal coverings of those root graded Lie
algebras. It turns out that the classification of those root
graded Lie algebras played a crucial role in classifying the newly
developed extended affine Lie algebras (see [BGKN] and [AG]). All
affine Kac-Moody Lie algebras except $A_{2l}^{(2)}$ are examples
of Lie algebras graded by reduced finite root systems.

To include the twisted affine Lie algebra $A_{2l}^{(2)}$ and for
the purpose of the classification of the extended affine Lie
algebras of non-reduced types, [ABG2] introduced Lie algebras
graded by the non-reduced root system $BC_{{}_N}$.
$BC_{{}_N}$-graded Lie algebras do appear not only in the extended
affine Lie algebras (see [AABGP]) including the twisted affine Lie
algebra $A_{2l}^{(2)}$ but also in the finite-dimensional
isotropic simple Lie algebras studied by Seligman [S]. The other
important examples include the ``odd symplectic" Lie algebras
studied by Gelfand-Zelevinsky [GeZ], Maliakas [Ma] and Proctor
[P].

       The Cliffold(or Weyl) algebras have natural representations on the exterior(or symmetric)
       algebras of polynomials over half of generators. Those representations are important in
       quantum and statistical mechanics where the generators are interpreted as operators which
       create or annihilate particles and satisfy Fermi(or Bose) statistics.
       Fermionic representations for the affine Kac-Moody Lie algebras were first
       developed by Frenkel \cite{F1} and Kac-Peterson \cite{KP} independently. Feingold-Frenkel \cite{FF}
       systematically constructed representations for all
classical affine
       Lie algebras by using Clifford or Weyl algebras with infinitely many generators.
       [G] constructed
       bosonic and fermionic representations for the extended affine Lie algebra
       $\widetilde{gl_N(\mathbb{C}_q)}$, where $\mathbb{C}_q$ is the quantum torus in two variables.
       Thereafter Lau [L] gave a more general construction.

       In this paper, we will construct fermions depending on the parameter $q$ which will lead to representations for
       some $BC_{{}_N}$-graded Lie algebras coordinatized by quantum tori with nontrivial central extensions. Since $C_{{}_N}$-graded Lie algebras are also $BC_{{}_N}$-graded
       Lie algebras we will treat bosons as well in a unified way.

       The organization of the paper is as follows. In Section 1, we review the definition of
       $BC_{{}_N}$-graded Lie algebras and give examples of $BC_{{}_N}$-graded Lie algebras
       which are
       subalgebras of $\widehat{gl_{2N}(\mathbb{C}_q)}$ or
       $\widehat{gl_{2N+1}(\mathbb{C}_q)}$.
        In Section 2, we use fermions or bosons to construct representations for those examples of $BC_{{}_N}$-graded Lie
       algebras by using Clifford or Weyl algebras with infinitely many
       generators. Although we get $BC_{{}_N}$-graded Lie algebras
       with the grading subalgebras of type $B_{{}_N}, C_{{}_N}$ and
       $D_{{}_N}$, there is only one which is a genuine
       $BC_{{}_N}$-graded Lie algebra arising from the fermionic
       construction.

       Throughout this paper, we denote the field of complex numbers and the ring of integers by
       $\mathbb{C}$ and $\mathbb{Z}$ respectively.

\section{BC$_{{}_{\text N}}$-graded Lie Algebras}

        In this section, we first recall the definition of quantum tori and
       $BC_{{}_N}$-graded Lie algebras. We then go on to provide some examples of BC$_N$-graded Lie algebras.
       For more information
       on BC$_N$-graded Lie algebras, see \cite{ABG2}.

       Let $q$ be a non-zero complex number. A quantum torus associated to $q$
       (see \cite{M}) is the unital associative $\mathbb{C}$-algebra $\mathbb{C}_q[x^{\pm},y^{\pm}]$
       (or simply $\mathbb{C}_q$) with generators $x^{\pm}$,$y^{\pm}$ and relations
           \begin{equation}
                xx^{-1}=x^{-1}x=yy^{-1}=y^{-1}y=1  \quad and \quad yx=qxy.
           \end{equation}
       Then
           \begin{equation}
                x^{m}y^{n}x^{p}y^{s}=q^{np}x^{m+p}y^{n+s}
           \end{equation}
       and
           \begin{equation}
                \mathbb{C}_q=\sum_{m,n \in \mathbb{Z}}\oplus\mathbb{C}x^my^n.
           \end{equation}
       Set $\Lambda(q)=\{n \in \mathbb{Z}|q^n=1\}$.
       From \cite{BGK} we see that $[\mathbb{C}_q,\mathbb{C}_q]$ has a basis consisting of monomials $x^my^n$ for
       $m \notin \Lambda(q)$ or $n \notin \Lambda(q) $.

       Let $\mbox{ }\bar{} \mbox{ }$ be the anti-involution on
       $\mathbb{C}_q$ given by
           \begin{equation}
           \bar{x}=x,\quad \bar{y}=y^{-1}.
           \end{equation}
       We have $\mathbb{C}_q=\mathbb{C}_q^+\oplus\mathbb{C}_q^-$, where
       $\mathbb{C}_q^{\pm}
               =\{s \in \mathbb{C}_q|\bar{s}=\pm s\}$, then \\
           \parbox{1cm}{\begin{eqnarray}\end{eqnarray}}\hfill \parbox{13.66cm}
               {\begin{eqnarray*}& \mathbb{C}_q^+=span\{x^my^n + \overline{x^my^n}|m \in \mathbb{Z},n \geq 0\}, & \\
                                  & \mathbb{C}_q^-=span\{x^my^n - \overline{x^my^n}|m \in \mathbb{Z},n>0\}. & \hfill
                \end{eqnarray*}}
       Now we form a central extension of $gl_r(\mathbb{C}_q)$ (cf. \cite{G}),
           \begin{equation}
               \widehat{gl_r(\mathbb{C}_q)}=gl_r(\mathbb{C}_q) \oplus
               \Bigl(\sum_{n \in \Lambda(q)}\oplus \mathbb{C}c(n) \Bigr) \oplus \mathbb{C}c_y
           \end{equation}
       with Lie bracket
           \begin{eqnarray}
                [e_{ij}(x^my^n),e_{kl}(x^py^s)]
                &=&\delta_{jk}q^{np}e_{il}(x^{m+p}y^{n+s})-\delta_{il}q^{ms}e_{kj}(x^{m+p}y^{n+s}) \nonumber  \\
                   &&+m q^{np}\delta_{jk}\delta_{il}\delta_{m+p,0}\delta_{\overline{n+s},\overline{0}}
                   c(n+s)  \\
                                  &&+n q^{np}\delta_{jk}\delta_{il}\delta_{m+p,0}\delta_{n+s,0}
                                  c_y \nonumber
           \end{eqnarray}
       for $m,p,n,s \in \mathbb{Z}$, where $c(u)$, for $u \in \Lambda(q)$ and $c_y$ are
       central elements of $\widehat{gl_r(\mathbb{C}_q)}$, $\overline{t} \mbox{ means }
       \overline{t} \in \mathbb{Z}/\Lambda(q),\mbox{for }t \in \mathbb{Z}$.

       Next we recall the definition of BC$_N$-graded Lie algebra and construct three types of
       BC$_N$-graded Lie algebras. Let       \\
           \parbox{1cm}{\begin{eqnarray}\end{eqnarray}}\hfill \parbox{13.66cm}
            {\begin{eqnarray*}
            &&\Delta_B=\{\pm \epsilon_i \pm \epsilon_j |1 \leq i \neq j \leq N\}\cup \{\pm \epsilon_i|i=1,\cdots ,N\} \\
            &&\Delta_C=\{\pm \epsilon_i \pm \epsilon_j |1 \leq i \neq j \leq N\}\cup \{\pm 2\epsilon_i|i=1,\cdots ,N\} \\
            &&\Delta_D=\{\pm \epsilon_i \pm \epsilon_j |1 \leq i \neq j \leq N\}. \hfill
                \end{eqnarray*}}
       be root systems of type $B$,$C$ and $D$ respectively, and
           \begin{equation}
              \Delta=\{\pm \epsilon_i \pm \epsilon_j |1 \leq i \neq j \leq N\}\cup
              \{\pm \epsilon_i,\pm 2\epsilon_i|i=1,\cdots ,N\}
           \end{equation}
       be a root system of type $BC_{{}_N}$ in the sense of Bourbaki [Bo, Chapitre VI].

           \begin{defini}
                 [BC$_N$-graded Lie Algebras]A Lie algebra L over a field $\mathbb{F}$ of characteristic 0 is
                 $\underline{\mbox{graded by the root system $BC_N$}}$ or is $\underline{\mbox{$BC_N$-graded}}$ if
              \begin{description}
                 \item[(i) ] L contained as a subalgebra a finite-dimentional split ``simple" Lie algebra
                 $\mathfrak{g}=\mathfrak{h}\oplus \bigoplus_{\mu \in \Delta_X}\mathfrak{g}_{\mu}$ whose root system
                 relative to a split Cartan subalgebra $\mathfrak{h}=\mathfrak{g}_0$ is $\Delta_X$, X=B,C, or D;
                 \item[(ii)] $L=\bigoplus_{\mu \in \Delta \cup \{0\}}L_{\mu}$, where $L_{\mu}=
                 \{x \in L|[h,x]=\mu(h)x, \mbox{ for all } h \in \mathfrak{h}\}$ for $\mu \in \Delta \cup \{0\}$,
                 and $\Delta$ is the root system $BC_N$ as in (1.9); and
                 \item[(iii)] $L_0=\sum_{\mu \in
                 \Delta}[L_{\mu},L_{-\mu}]$.
              \end{description}

            \end{defini}

            In Definition 1.1 the word simple is in quotes, because
in every case but two the Lie algebra $\frak g$ associated with
$\Delta_{\text{X}}$ is simple; the sole exceptions being when
$\Delta_{\text{X}}$ = D$_2$ or D$_1$.  The D$_2$ root system is the
same as $\text{A}_1 \times \text{A}_1$, and $\frak g$ is the sum
$\frak{g} = \frak{g}^{(1)} \oplus \frak{g}^{(2)}$ of two copies of
$\frak{sl}_2$ in this case. In the D$_1$ case, $\frak{g} =
\mathbb{F} h$, a one-dimensional subalgebra.

We refer to $\frak g$ as the $\underline {\text{grading
subalgebra}}$ of $L$, and we say $L$ is $\underline
{\text{$BC_{{}_N}$-graded}}$ \ $\underline {\text{with grading
subalgebra $\frak g$ of type X$_{{}_N}$}} $ (where X = B, C, or D)
to mean that the root system of $\frak g$ is of type X$_{{}_N}$.

 Any Lie algebra which is graded by a finite
root system of type B$_{{}_N}$, C$_{{}_N}$, or D$_{{}_N}$ is also
BC$_{{}_N}$-graded with grading subalgebra of type B$_{{}_N}$,
C$_{{}_N}$, or D$_{{}_N}$ respectively. For such a Lie algebra $L$,
the space $L_\mu = (0)$ for all $\mu$ not in $\Delta_B$, $\Delta_C$,
or $\Delta_D$ respectively.

\subsection{Type C and D }

 For BC$_N$-graded Lie algebras with grading subalgebra of
       type $C_N(\rho=-1)$ and $D_N(\rho=1)$, we put
            \begin{equation*}G=\left(
                 \begin{array}{cc}
                      0 & I_N \\
                      \rho I_N & 0 \\
                 \end{array}
                      \right) \in M_{2N}(\mathbb{C}_q).
            \end{equation*}
       Then, $G$ is an invertible $2N \times 2N-$matrix and $\bar{G}^t=\rho G$. Using the
       matrix $G$, we define a map
       $$^*:M_{2N}(\mathbb{C}_q) \rightarrow
       M_{2N}(\mathbb{C}_q)\text{
       by } A^*=G^{-1}\bar{A}^tG.$$
        Since $\bar{G}^t=\rho G$, $^*$ is an
       involution of the associative algebra
       $M_{2N}(\mathbb{C}_q)$. As in \cite{AABGP}, we define
            \begin{equation*}
                 S_{\rho}(M_{2N}(\mathbb{C}_q),^*)=\{A \in M_{2N}(\mathbb{C}_q):A^*=-A \}
            \end{equation*}
       in which case $S_{\rho}(M_{2N}(\mathbb{C}_q),^*)$ is a Lie
       subalgebra of $gl_{2N}(\mathbb{C}_q)$ over $\mathbb{C}$. The
       general form of a matrix in $S_{\rho}(M_{2N}(\mathbb{C}_q),^*)$ is
            \begin{equation}     \left(
                 \begin{array}{cc}
                      A & S \\
                      T & -\bar{A}^t\\
                 \end{array}
                      \right) \quad \mbox{with  } \bar{S}^t=-\rho S \quad \mbox{and}
                      \quad     \bar{T}^t=-\rho T
            \end{equation}
       where $A,S,T \in M_N(\mathbb{C}_q)$. Then the Lie
       algebra
            \begin{equation*}
                 \mathcal{G}_{\rho}=[S_{\rho}(M_{2N}(\mathbb{C}_q),^*),S_{\rho}(M_{2N}(\mathbb{C}_q),^*)]
            \end{equation*}
            is a BC$_{N}$-graded Lie algebra with grading
            subalgebra of type $C_N(\rho=-1)$ and $D_N(\rho=1)$.
       Using the method in \cite{AABGP}, we easily know that
$$\mathcal{G}_{\rho}=\bigl\{Y \in
S_{\rho}(M_{2N}(\mathbb{C}_q),^*)|tr(Y)
                 \equiv 0\ \mbox{mod }
                 [\mathbb{C}_q,\mathbb{C}_q]\bigr\}.$$

       We put
            \begin{equation}
                 \mathcal{H}=
                 \Bigl\{\sum_{i=1}^Na_i(e_{ii}-e_{N+i,N+i})|a_i \in
                 \mathbb{C}\Bigr\},
            \end{equation}
       then $\mathcal{H}$ is a $N$-dimensional abelian subalgebra of $\mathcal{G}_{\rho}$. Defining $\epsilon_i \in
       \mathcal{H}^*,i=1,\cdots,N$, by
            \begin{equation}
                 \epsilon_i\biggl(\sum_{j=1}^Na_j(e_{jj}-e_{N+j,N+j}) \biggr)=a_i
            \end{equation}
       for $i=1,\cdots,N.$ Putting $\mathcal{G}_{\alpha}=\{x\in \mathcal{G}_{\rho}|[h,x]=\alpha(h)x,
       \mbox{ for all } h \in \mathcal{H}\}$ as usual, we have
            \begin{equation}
                 \mathcal{G}_{\rho}=\mathcal{G}_0 \oplus \sum_{i \neq j} \mathcal{G}_{\epsilon_i-\epsilon_j}
                 \oplus \sum_{i<j}(\mathcal{G}_{\epsilon_i+\epsilon_j} \oplus
                 \mathcal{G}_{-\epsilon_i-\epsilon_j})\oplus \sum_{i}
                 (\mathcal{G}_{2\epsilon_i} \oplus \mathcal{G}_{-2\epsilon_i})
            \end{equation}
       where \\
            \parbox{1cm}{\begin{eqnarray}\end{eqnarray}}\hfill \parbox{13.66cm}
               {\begin{eqnarray*}
               &\mathcal{G}_{\epsilon_i-\epsilon_j}=span_{\mathbb{C}}\{\tilde{f}_{ij}(m,n)=x^my^ne_{ij}-
               \overline{x^my^n}e_{N+j,N+i} |m,n \in \mathbb{Z}\}, & \\
               &\mathcal{G}_{\epsilon_i+\epsilon_j}=span_{\mathbb{C}}\{ \tilde{g}_{ij}(m,n)=x^my^ne_{i,N+j}-
               \rho \overline{x^my^n}e_{j,N+i}|m,n \in \mathbb{Z}\},& \\
               &\mathcal{G}_{-\epsilon_i-\epsilon_j}=span_{\mathbb{C}}\{ \tilde{h}_{ij}(m,n)=\rho x^my^ne_{N+i,j}-
               \overline{x^my^n}e_{N+j,i} |m,n \in \mathbb{Z}\}, & \\
               &\mathcal{G}_{2\epsilon_i}=span_{\mathbb{C}}\{\tilde{g}_{ii}(m,n)=
               (x^my^n-\rho \overline{x^my^n})e_{i,N+i}|m,n \in \mathbb{Z} \}, & \\
               &\mathcal{G}_{-2\epsilon_i}=span_{\mathbb{C}}\{\tilde{h}_{ii}(m,n)=
               (\rho x^my^n- \overline{x^my^n})e_{N+i,i}|m,n \in \mathbb{Z} \}, &\hfill
                \end{eqnarray*}}
       and
            \begin{equation*}
                 \mathcal{G}_0=span_{\mathbb{C}}\{\tilde{f}_{ii}(m,n)-\tilde{f}_{11}(m,n),\tilde{f}_{11}(p,s)
                 |2 \leq i \leq N,m,n \in \mathbb{Z},p \notin \Lambda(q)
               \mbox{ or } s \notin \Lambda(q) \}.
            \end{equation*}
       Note that $\tilde{g}_{ij}(m,n)=-\rho q^{-mn}\tilde{g}_{ji}(m,-n)$, $\tilde{h}_{ij}(m,n)=-\rho
       q^{-mn}\tilde{h}_{ji}(m,-n).$

       Now we form a central extension of $\mathcal{G}_{\rho}$
            \begin{equation}
                 \widehat{\mathcal{G}_{\rho}}=\mathcal{G}_{\rho} \oplus
                 \Bigl(\sum_{n \in \Lambda(q)}\oplus \mathbb{C}c(n) \Bigr) \oplus \mathbb{C}c_y
            \end{equation}
       with Lie brackets as (1.7).

       We have
            \begin{prop}
                 \begin{equation}
                      [\tilde{g}_{ij}(m,n),\tilde{g}_{kl}(p,s)]=0
                 \end{equation}
                 \begin{equation}
                      [\tilde{g}_{ij}(m,n),\tilde{f}_{kl}(p,s)]=-\delta_{il}q^{ms}\tilde{g}_{kj}(m+p,n+s)
                      +\rho \delta_{jl}q^{(s-n)m}\tilde{g}_{ki}(m+p,s-n)
                 \end{equation}
            \begin{eqnarray}
                 &&[\tilde{g}_{ij}(m,n),\tilde{h}_{kl}(p,s)] \nonumber \\
                 &=&-\delta_{ik}q^{-n(m+p)}\tilde{f}_{jl}(m+p,s-n)
                 +\rho \delta_{jk}q^{np}\tilde{f}_{il}(m+p,n+s)  \nonumber  \\
                 &&+\rho \delta_{il}q^{-(mn+np+ps)}\tilde{f}_{jk}(m+p,-(n+s))
                 -\delta_{jl}q^{(n-s)p}\tilde{f}_{ik}(m+p,n-s)      \\
                 &&+m\rho q^{np}\delta_{jk}\delta_{il}\delta_{m+p,0}\delta_{\overline{n+s},\overline{0}}(c(n+s)+c(-n-s))
                         \nonumber \\
                 &&-m\delta_{ik}\delta_{jl}\delta_{m+p,0}\delta_{\overline{n-s},\overline{0}}(c(n-s)+c(s-n)) \nonumber
            \end{eqnarray}
            \begin{eqnarray}
                 [\tilde{f}_{ij}(m,n),\tilde{f}_{kl}(p,s)]&=&\delta_{jk}q^{np}
                 \tilde{f}_{il}(m+p,n+s)-\delta_{il}q^{sm}\tilde{f}_{kj}(m+p,n+s) \nonumber   \\
                 &&+2mq^{np}\delta_{jk}\delta_{il}\delta_{m+p,0} \delta_{\overline{n+s},\overline{0}}c(n+s)
            \end{eqnarray}
            \begin{equation}
                 [\tilde{f}_{ij}(m,n),\tilde{h}_{kl}(p,s)]=
                 -\delta_{ik}q^{-n(m+p)}\tilde{h}_{jl}(m+p,s-n)-\delta_{il}q^{ms}\tilde{h}_{kj}(m+p,n+s)
            \end{equation}
            \begin{equation}
                 [\tilde{h}_{ij}(m,n),\tilde{h}_{kl}(p,s)]=0
            \end{equation}
       for all $m,p,n,s \in \mathbb{Z}$ and $1 \leq i,j,k,l \leq N$.
            \end{prop}
       \noindent\textbf{Proof. } We only check (1.18).
            \begin{eqnarray*}
                 &&[\tilde{g}_{ij}(m,n),\tilde{h}_{kl}(p,s)]\\
                 &=&[x^my^ne_{i,N+j}-\rho\overline{x^my^n}e_{j,N+i},\rho x^py^se_{N+k,l}-\overline{x^py^s}e_{N+l,k}] \\
                 &=&\rho[x^my^ne_{i,N+j},x^py^se_{N+k,l}]-[x^my^ne_{i,N+j},\overline{x^py^s}e_{N+l,k}]
                    -[\overline{x^my^n}e_{j,N+i},x^py^se_{N+k,l}]                    \\
                    &&+\rho[\overline{x^my^n}e_{j,N+i},\overline{x^py^s}e_{N+l,k}]    \\
                 &=&\rho\bigl(\delta_{jk}x^my^nx^py^se_{il}-\delta_{il}x^py^sx^my^ne_{N+k,N+j}
                     +mq^{np}\delta_{jk}\delta_{il}\delta_{m+p,0}\delta_{\overline{n+s},\overline{0}}c(n+s)\bigr)  \\
                   && -\bigl(\delta_{jl}x^my^n \overline{x^py^s}e_{ik}-\delta_{ki}\overline{x^py^s}x^my^n e_{N+l,N+j}
                      +m\delta_{jl}\delta_{ik}\delta_{m+p,0}\delta_{\overline{n-s},\overline{0}}c(n-s)\bigr)   \\
                 &&-\bigl(\delta_{ik}\overline{x^my^n}x^py^se_{jl}-\delta_{lj}x^py^s\overline{x^my^n}e_{N+k,N+i}
                    +m\delta_{jl}\delta_{ik}\delta_{m+p,0}\delta_{\overline{n-s},\overline{0}}c(s-n)\bigr)  \\
                 &&+\rho\bigl(\delta_{il}\overline{x^py^sx^my^n}e_jk-\delta_{kj}
                 \overline{x^my^nx^py^s}e_{N+l,N+i}
                 +mq^{np}\delta_{jl}\delta_{ik}\delta_{m+p,0}\delta_{\overline{n+s},\overline{0}}c(-n-s)\bigr)   \\
                 &&+\rho n\delta_{jk}\delta_{il}\delta_{m+p,0}\delta_{n+s,0}c_y
                 -n\delta_{jl}\delta_{ik}\delta_{m+p,0}\delta_{n-s,0}c_y
                 +n\delta_{jl}\delta_{ik}\delta_{m+p,0}\delta_{n-s,0}c_y  \\
                 &&-\rho n \delta_{jk}\delta_{il}\delta_{m+p,0}\delta_{n+s,0}c_y      \\
                 &=&-\delta_{ik}q^{-n(m+p)}\tilde{f}_{jl}(m+p,s-n)
                 +\rho \delta_{jk}q^{np}\tilde{f}_{il}(m+p,n+s)    \\
                 &&+\rho \delta_{il}q^{-(mn+np+ps)}\tilde{f}_{jk}(m+p,-(n+s))
                 -\delta_{jl}q^{(n-s)p}\tilde{f}_{ik}(m+p,n-s)      \\
                 &&+m\rho q^{np}\delta_{jk}\delta_{il}\delta_{m+p,0}\delta_{\overline{n+s},\overline{0}}(c(n+s)+c(-n-s))
                        \\
                 &&-m\delta_{ik}\delta_{jl}\delta_{m+p,0}\delta_{\overline{n-s},\overline{0}}(c(n-s)+c(s-n)).
            \end{eqnarray*}
        The proof of the others is similar.       $\hfill \blacksquare$

\subsection{Type B}

        For type $B$, we put
            \begin{equation*}G=\left(
                 \begin{array}{ccc}
                      1 & 0 & 0 \\
                      0 & 0 & I_N \\
                      0 & I_N & 0 \\
                 \end{array}
                      \right) \in M_{2N+1}(\mathbb{C}_q).
            \end{equation*}
       Then, $G$ is an invertible {\small $ (2N+1) \times (2N+1)$}-matrix and $\bar{G}^t=G$. Using the
       matrix $G$, we define a map
       $$^*:M_{2N+1}(\mathbb{C}_q) \rightarrow
       M_{2N+1}(\mathbb{C}_q)\text{
       by }A^*=G^{-1}\bar{A}^tG.$$
       Since $\bar{G}^t=G$, \mbox{ $^*$} is an
       involution of the associative algebra $M_{2N+1}(\mathbb{C}_q)$. As in \cite{AABGP}, we define
            \begin{equation*}
                 S(M_{2N+1}(\mathbb{C}_q),^*)=\{A \in M_{2N+1}(\mathbb{C}_q):A^*=-A \}
            \end{equation*}
       in which case $S(M_{2N+1}(\mathbb{C}_q),^*)$ is a Lie subalgebra of $gl_{2N+1}(\mathbb{C}_q)$
       over $\mathbb{C}$. The general form of a matrix in $S(M_{2N+1}(\mathbb{C}_q),^*)$ is
            \begin{equation}
                 \left(
            \begin{array}{ccc}
                 a & b_1 & b_2 \\
                 -\bar{b_2}^t  & A & S \\
                 -\bar{b_1}^t & T & -\bar{A}^t\\
            \end{array}
                 \right) \quad \mbox{with  }\bar{a}=-a \quad \bar{S}^t=-S \quad \mbox{and}
                 \quad     \bar{T}^t=-T
            \end{equation}
       where $A,S,T \in M_N(\mathbb{C}_q)$. Then the Lie
       algebra
            \begin{equation*}
                 \mathcal{G'}=[S(M_{2N+1}(\mathbb{C}_q),^*),S(M_{2N+1}(\mathbb{C}_q),^*)]
            \end{equation*}
            is a BC$_N$-graded Lie algebra with grading subalgebra
            of type B$_N$.
       Following from \cite{AABGP}, we easily know that
$$\mathcal{G'}=\{Y \in S(M_{2N+1}(\mathbb{C}_q),^*)|tr(Y)
                 \equiv 0\ \mbox{mod }
                 [\mathbb{C}_q,\mathbb{C}_q]\}$$
       As in Section 1.1, we set
            \begin{equation}
                 \mathcal{H}'=\Bigl\{\sum_{i=1}^Na_i(e_{ii}-e_{N+i,N+i})|a_i \in
                 \mathbb{C}\Bigr\},
            \end{equation}
      then $\mathcal{H}'$ is a $N$-dimensional abelian subalgebra of $\mathcal{G}'$. Defining $\epsilon_i \in
       \mathcal{H'}^*,i=1,\cdots,N$, by
            \begin{equation}
                 \epsilon_i\biggl(\sum_{j=1}^la_j(e_{jj}-e_{N+j,N+j}) \biggr)=a_i
            \end{equation}
       for $i=1,\cdots,N.$ Putting $\mathcal{G}'_{\alpha}=\{x\in \mathcal{G}'|[h,x]=\alpha(h)x,
       \mbox{ for all } h \in \mathcal{H'}\}$ as usual, we have
                   \begin{equation}
                 \mathcal{G}'=\mathcal{G}'_0 \oplus \sum_{i \neq j} \mathcal{G}'_{\epsilon_i-\epsilon_j}
                 \oplus \sum_{i<j}(\mathcal{G}'_{\epsilon_i+\epsilon_j} \oplus \mathcal{G}'_{-\epsilon_i-\epsilon_j})
                 \oplus \sum_i (\mathcal{G}'_{\epsilon_i} \oplus \mathcal{G}'_{-\epsilon_i}
                 \oplus \mathcal{G}'_{2\epsilon_i} \oplus \mathcal{G}'_{-2\epsilon_i})
            \end{equation}
       where \\
            \parbox{1cm}{\begin{eqnarray}\end{eqnarray}}\hfill \parbox{13.66cm}
               {\begin{eqnarray*}
               &\mathcal{G}'_{\epsilon_i-\epsilon_j}=span_{\mathbb{C}}\{\tilde{f}_{ij}(m,n)=x^my^ne_{ij}-
               \overline{x^my^n}e_{N+j,N+i} |m,n \in \mathbb{Z}\}, & \\
               &\mathcal{G}'_{\epsilon_i+\epsilon_j}=span_{\mathbb{C}}\{ \tilde{g}_{ij}(m,n)=x^my^ne_{i,N+j}-
               \overline{x^my^n}e_{j,N+i}|m,n \in \mathbb{Z} \},& \\
               &\mathcal{G}'_{-\epsilon_i-\epsilon_j}=span_{\mathbb{C}}\{ \tilde{h}_{ij}(m,n)=x^my^ne_{N+i,j}-
               \overline{x^my^n}e_{N+j,i} |m,n \in \mathbb{Z} \}, & \\
               &\mathcal{G}'_{2\epsilon_i}=span_{\mathbb{C}}\{\tilde{g}_{ii}(m,n)=
               (x^my^n-\overline{x^my^n})e_{i,N+i}|m,n \in \mathbb{Z} \}, & \\
               &\mathcal{G}'_{-2\epsilon_i}=span_{\mathbb{C}}\{\tilde{h}_{ii}(m,n)=
               (x^my^n- \overline{x^my^n})e_{N+i,i}|m,n \in \mathbb{Z} \},&\\
               &\mathcal{G}'_{\epsilon_i}=span_{\mathbb{C}}\{\tilde{e}_i(m,n)=x^my^ne_{i,0}-
               \overline{x^my^n}e_{0,N+i}|m,n \in \mathbb{Z} \}, & \\
               &\mathcal{G}'_{-\epsilon_i}=span_{\mathbb{C}}\{\tilde{e}_i^*(m,n)=x^my^ne_{N+i,0}-
               \overline{x^my^n}e_{0,i}|m,n \in \mathbb{Z}\}, &\hfill
                \end{eqnarray*}}
       and
            \begin{equation*}
                 \mathcal{G}'_0=span_{\mathbb{C}}\{\tilde{f}_{ii}(m,n)-\tilde{e}_0(m,n),\tilde{e}_0(p,s)
                 |1 \leq i \leq N,m,n \in \mathbb{Z},p \notin \Lambda(q)
               \mbox{ or } s \notin \Lambda(q)\},
            \end{equation*}
       where $\tilde{e}_0(m,n)=(x^my^n-\overline{x^my^n})e_{0,0}$.

       Next we form a central extension of $\mathcal{G}'$
            \begin{equation}
                 \widehat{\mathcal{G}'}=\mathcal{G}' \oplus
                 \Bigl(\sum_{n \in \Lambda(q)}\oplus \mathbb{C}c(n) \Bigr) \oplus \mathbb{C}c_y
            \end{equation}
       with Lie brackets as (1.7).
            \begin{rem}
                Note that the index of the matrices in $M_{2N+1}(\mathbb{C}_q)$
                 ranges from $0$ to $2N$.
            \end{rem}
       Now we have
            \begin{prop}
                 \begin{equation}
                      [\tilde{g}_{ij}(m,n),\tilde{g}_{kl}(p,s)]=0
                 \end{equation}
                 \begin{equation}
                      [\tilde{g}_{ij}(m,n),\tilde{f}_{kl}(p,s)]=-\delta_{il}q^{ms}\tilde{g}_{kj}(m+p,n+s)
                      +\delta_{jl}q^{(s-n)m}\tilde{g}_{ki}(m+p,s-n)
                 \end{equation}
            \begin{eqnarray}
                 &&[\tilde{g}_{ij}(m,n),\tilde{h}_{kl}(p,s)]\\
                 &=&-\delta_{ik}q^{-n(m+p)}\tilde{f}_{jl}(m+p,s-n)
                 +\delta_{jk}q^{np}\tilde{f}_{il}(m+p,n+s)  \nonumber  \\
                 &&+\delta_{il}q^{-(mn+np+ps)}\tilde{f}_{jk}(m+p,-(n+s))
                 -\delta_{jl}q^{(n-s)p}\tilde{f}_{ik}(m+p,n-s)   \nonumber   \\
                 &&+m q^{np}\delta_{jk}\delta_{il}\delta_{m+p,0}\delta_{\overline{n+s},\overline{0}}(c(n+s)+c(-n-s))
                         \nonumber \\
                 &&-m\delta_{ik}\delta_{jl}\delta_{m+p,0}\delta_{\overline{n-s},\overline{0}}(c(n-s)+c(s-n)) \nonumber
            \end{eqnarray}
            \begin{equation}
                 [\tilde{g}_{ij}(m,n),\tilde{e}_{k}(p,s)]=0
            \end{equation}
            \begin{equation}
                 [\tilde{g}_{ij}(m,n),\tilde{e}_k^*(p,s)]=-\delta_{ik}q^{-n(m+p)}
                 \tilde{e}_j(m+p,s-n)+ \delta_{jk}q^{np}\tilde{e}_i(m+p,n+s)
            \end{equation}
            \begin{equation}
                 [\tilde{g}_{ij}(m,n),\tilde{e}_{0}(p,s)]=0
            \end{equation}
            \begin{eqnarray}
                 [\tilde{f}_{ij}(m,n),\tilde{f}_{kl}(p,s)]&=&\delta_{jk}q^{np}
                 \tilde{f}_{il}(m+p,n+s)-\delta_{il}q^{sm}\tilde{f}_{kj}(m+p,n+s) \nonumber  \\
                  &&+2mq^{np}\delta_{jk}\delta_{il}\delta_{m+p,0}\delta_{\overline{n+s},\overline{0}}c(n+s)
            \end{eqnarray}
            \begin{equation}
                 [\tilde{f}_{ij}(m,n),\tilde{h}_{kl}(p,s)]=-\delta_{ik}q^{-n(m+p)}
                 \tilde{h}_{jl}(m+p,s-n)-\delta_{il}q^{ms}\tilde{h}_{kj}(m+p,n+s)
            \end{equation}
            \begin{equation}
                 [\tilde{f}_{ij}(m,n),\tilde{e}_k(p,s)]=\delta_{jk}q^{np}\tilde{e}_i(m+p,n+s)
            \end{equation}
            \begin{eqnarray}
                 &[\tilde{f}_{ij}(m,n),\tilde{e}_k^*(p,s)]=-\delta_{ik}q^{-n(m+p)}\tilde{e}_j^*(m+p,s-n)&
            \end{eqnarray}
            \begin{equation}
                 [\tilde{f}_{ij}(m,n),\tilde{e}_{0}(p,s)]=0
            \end{equation}
            \begin{equation}
                 [\tilde{h}_{ij}(m,n),\tilde{h}_{kl}(p,s)]=0
            \end{equation}
            \begin{equation}
                 [\tilde{h}_{ij}(m,n),\tilde{e}_k(p,s)]=\delta_{jk}q^{np}\tilde{e}_i^*(m+p,n+s)-
                   \delta_{ik}q^{-n(m+p)}\tilde{e}_j^*(m+p,s-n)
            \end{equation}
            \begin{equation}
                 [\tilde{h}_{ij}(m,n),\tilde{e}^*_{k}(p,s)]=0
            \end{equation}
            \begin{equation}
                 [\tilde{h}_{ij}(m,n),\tilde{e}_{0}(p,s)]=0
            \end{equation}
            \begin{equation}
                 [\tilde{e}_i(m,n),\tilde{e}_k(p,s)]=q^{m(s-n)}\tilde{g}_{ki}(m+p,s-n)
            \end{equation}
            \begin{eqnarray}
                 &&[\tilde{e}_i(m,n),\tilde{e}_k^*(p,s)] \nonumber  \\
                 &=&-\delta_{ik}q^{-n(m+p)}\tilde{e}_0(m+p,s-n)-q^{p(n-s)}\tilde{f}_{ik}(m+p,n-s) \\
                 &&+m\delta_{ik}\delta_{m+p,0}\delta_{\overline{n-s},\overline{0}}(c(n-s)+c(s-n)) \nonumber
            \end{eqnarray}
            \begin{equation}
                 [\tilde{e}_i(m,n),\tilde{e}_0(p,s)]=q^{np}\tilde{e}_i(m+p,n+s)-q^{p(n-s)}\tilde{e}_i(m+p,n-s)
            \end{equation}
            \begin{equation}
                 [\tilde{e}_i^*(m,n),\tilde{e}_k^*(p,s)]=q^{m(s-n)}\tilde{h}_{ki}(m+p,s-n)
            \end{equation}
            \begin{equation}
                 [\tilde{e}_i^*(m,n),\tilde{e}_0(p,s)]=q^{np}\tilde{e}_i^*(m+p,n+s)-q^{p(n-s)}\tilde{e}_i^*(m+p,n-s)
            \end{equation}
            \begin{eqnarray}
                 &&[\tilde{e}_0(m,n),\tilde{e}_0(p,s)]  \\
                 &=&(q^{np}-q^{sm})\tilde{e}_0(m+p,n+s)+(q^{m(s-n)}-q^{-n(m+p)})\tilde{e}_0(m+p,s-n) \nonumber \\
                 &&+mq^{np}\delta_{m+p,0}\delta_{\overline{n+s},\overline{0}}(c(n+s)+c(-n-s))\nonumber  \\
                 &&-m\delta_{m+p,0}\delta_{\overline{n-s},\overline{0}}(c(n-s)+c(s-n)) \nonumber
            \end{eqnarray}
       for all $m,p,n,s \in \mathbb{Z}$ and $1 \leq i,j,k,l \leq N$.
            \end{prop}

                 The proof of Proposition 1.2 is similar to Proposition 1.1.

\begin{rem} Note that unlike (1.4),  the anti-involution in
[AABGP] is given by
$$\bar{x}=\pm x, \ \ \bar{y} =\pm y.$$
   \end{rem}

\section{Representations}

        In this section, we follow the idea in \cite{G} and \cite{FF}  to construct representations
       for the three types of BC$_N$-graded Lie algebras which are given in Section 1.

       Let $\mathcal{R}$ be an associative algebra. Let $\rho=\pm 1$. We define a $\rho$-bracket on
       $\mathcal{R}$ as follow:
            \begin{equation}
                 \{a,b\}_{\rho}=ab+\rho ba,\quad a,b \in \mathcal{R}.
            \end{equation}
       It is easy to see that
            \begin{equation}
                 \{a,b\}_{\rho}=\rho \{b,a\}_{\rho} \mbox{  and } [ab,c]=a\{b,c\}_{\rho}-\rho
                 \{a,c\}_{\rho}b
            \end{equation}
       for $a,b,c \in \mathcal{R}$, where $[a,b]=\{a,b\}_{-1}$ is the Lie bracket.

\subsection{Type C and D}

       Define $\mathfrak{a}$ to be the unital associative algebra with $2N$ generators
       $a_i,a_i^*,1 \leq i \leq N $, subject to relations
            \begin{equation}
                 \{a_i,a_j\}_{\rho}=\{a_i^*,a_j^*\}_{\rho}=0, \quad \mbox{and}
                                    \quad \{a_i,a_j^*\}_{\rho}=\rho \delta_{ij}.
            \end{equation}
       Let the associative algebra $\alpha (N,\rho)$ be generated by
            \begin{equation}
                 \{u(m)|u \in \bigoplus_{i=1}^N (\mathbb{C}a_i \oplus \mathbb{C}a^*_i),m \in \mathbb{Z}\}
            \end{equation}
       with the relations
            \begin{equation}
                 \{u(m),v(n)\}_{\rho}=\{u,v\}_{\rho}\delta_{m+n,0}.
            \end{equation}
       We now define the normal ordering as in \cite{FF}(see also \cite{F2}). \\
            \parbox{1cm}{\begin{eqnarray}\end{eqnarray}}\hfill \parbox{13.66cm}
               {\begin{eqnarray*} :u(m)v(n):&=&\left\{
                   \begin{array}
                        {r@{ }l}  & u(m)v(n) \hspace{3.45cm}                             \mbox{if } n>m,
                               \\ & \frac{1}{2}\bigl(u(m)v(n)-\rho v(n)u(m)\bigr)  \quad \mbox{if } m=n,
                               \\ & -\rho v(n)u(m) \hspace{2.9cm}                        \mbox{if } m>n,
                   \end{array} \right.  \\&=&-\rho:v(n)u(m): \hfill
                \end{eqnarray*}}
       for $n,m \in \mathbb{Z},u,v \in \mathfrak{a}.$ Set
            \begin{equation}
                 \theta(n)=\left\{
                   \begin{array}
                        {r@{ \quad }l}
                         1, & \mbox{for } n>0, \\
                         \frac{1}{2} ,& \mbox{for } n=0,  \\
                         0, & \mbox{for } n<0,
                   \end{array} \right. \mbox{   then  }     1-\theta(n)=\theta(-n).
            \end{equation}
       We have \\
            \parbox{1cm}{\begin{eqnarray}\end{eqnarray}}\hfill \parbox{13.66cm}
               {\begin{eqnarray*}& :a_i(m)a_j(n):=a_i(m)a_j(n)=-\rho a_j(n)a_i(m),& \\
                                  & :a_i^*(m)a_j^*(n):=a_i^*(m)a_j^*(n)=-\rho a_j^*(n)a_i^*(m).& \hfill
                \end{eqnarray*}}
       and \\
            \parbox{1cm}{\begin{eqnarray}\end{eqnarray}}\hfill \parbox{13.66cm}
               {\begin{eqnarray*}
                     &a_i(m)a_j^*(n)=:a_i(m)a_j^*(n):+\rho \delta_{ij}\delta_{m+n,0}   \theta(m-n), & \\
                     & a_j^*(n)a_i(m)=-\rho:a_i(m)a_j^*(n):+\delta_{ij}\delta_{m+n,0}\theta(n-m).&
                                  \hfill
                \end{eqnarray*}}

       It follows from (2.2) that  \\
            \parbox{1cm}{\begin{eqnarray}\end{eqnarray}}\hfill \parbox{13.66cm}
               {\begin{eqnarray*}&&[a_i(m)a_j(n),a_k(p)]=0,  \\
                                 &&[a_i(m)a_j(n),a_k^*(p)]=-\delta_{ik} \delta_{m+p,0}a_j(n)+\rho \delta_{jk}
                                 \delta_{n+p,0} a_i(m), \\
                                 && [a_i(m)a_j^*(n),a_k(p)]=\delta_{jk} \delta_{n+p,0}a_i(m),  \\
                                 && [a_i(m)a_j^*(n),a_k^*(p)]=-\delta_{ik}\delta_{m+p,0}a_j^*(n),  \\
                                 &&[a_i^*(m)a_j^*(n),a_k(p)]=
                                     \delta_{jk} \delta_{n+p,0}a_i^*(m)-\rho \delta_{ik} \delta_{m+p,0} a_j^*(n), \\
                                 &&[a_i^*(m)a_j^*(n),a_k^*(p)]=0,\\   \hfill
                \end{eqnarray*}}
       for $m,n,p \in \mathbb{Z},1 \leq i,j,k \leq N.$

       Let $\alpha(N,\rho)^+$ be the subalgebra generated by $a_i(n),a_j^*(m),a_k^*(0)$,
       for $n,m> \nolinebreak 0$,
       and $1 \leq i,j,k \leq N.$ Let $\alpha(N,\rho)^-$ be the subalgebra generated by
       $a_i(n)$, $a_j^*(m)$, $a_k(0),$
       for $n,m<0,$ and $1 \leq i,j,k \leq N.$ Those generators in $\alpha(N,\rho)^+$ are called annihilation
       operators while those in $\alpha(N,\rho)^-$ are called creation operators. Let $V(N,\rho)$ be a simple
       $\alpha(N,\rho)$-module containing an element $v_0$, called a``vacuum vector", and satisfying
            \begin{equation}
                 \alpha(N,\rho)^+v_0=0.
            \end{equation}
       So all annihilation operators kill $v_0$ and
            \begin{equation}
                 V(N,\rho)=\alpha(N,\rho)^-v_0.
            \end{equation}

       Now we are in the position to construct a class of fermions (if $\rho=1$) or bosons (if $\rho=-1$) on $V(N,\rho)$. For
       any $m,n \in \mathbb{Z},1 \leq i,j \leq N,$ set
            \begin{equation}
                 f_{ij}(m,n)=\sum_{s \in \mathbb{Z}} q^{-ns}:a_i(m-s)a_j^*(s):,
            \end{equation}
            \begin{equation}
                 g_{ij}(m,n)=\sum_{s \in \mathbb{Z}} q^{-ns}:a_i(m-s)a_j(s):,
            \end{equation}
            \begin{equation}
                 h_{ij}(m,n)=\sum_{s \in \mathbb{Z}} q^{-ns}:a_i^*(m-s)a_j^*(s):.
            \end{equation}
       Although $f_{ij}(m,n),g_{ij}(m,n) \mbox{ and } h_{ij}(m,n)$ are infinite sums, they are
       well-defined as operators on $V(N,\rho)$. Indeed, for any vector $v \in V(N,\rho)=\alpha(N,\rho)^-v_0,$
       only finitely many terms in (2.13)-(2.15) can make a non-zero contribution to
       $g_{ij}(m,n)v$, $f_{ij}(m,n)v$ and $h_{ij}(m,n)v$.

            \begin{lemma}
                 We have \\
                 \parbox{1cm}{\begin{eqnarray}\end{eqnarray}}\hfill \parbox{13.66cm}
                 {\begin{eqnarray*}& g_{ij}(m,n)=-\rho q^{-mn}g_{ji}(m,-n), & \\
                                  & h_{ij}(m,n)=-\rho q^{-mn}h_{ji}(m,-n). & \hfill
                  \end{eqnarray*}}
       for m,n,p,s $\in \mathbb{Z}$ and $1 \leq i,j,k,l \leq N$.
            \end{lemma}
       \noindent\textbf{Proof. } We only prove for $g_{ij}(m,n)$. The proof of $h_{ij}(m,n)$ is similar.
            \begin{eqnarray*}
                 g_{ij}(m,n)&=&\sum_{s \in \mathbb{Z}}q^{-ns}:a_i(m-s)a_j(s):                \\
                 &=&-\rho \sum_{s \in \mathbb{Z}}q^{-ns}:a_j(s)a_i(m-s):                     \\
                 &=&-\rho \sum_{s \in \mathbb{Z}}q^{-n(m-s)} :a_j(m-s)a_i(s):                \\
                 &=&-\rho q^{-mn}g_{ji}(m,-n).
            \end{eqnarray*} $\hfill \blacksquare$

            \begin{lemma} We have
                 \begin{equation}
                      [g_{ij}(m,n),a_k(p)]=0,
                 \end{equation}
                 \begin{equation}
                      [g_{ij}(m,n),a_k^*(p)]=-\delta_{ik}q^{-n(m+p)}a_j(m+p)+\rho \delta_{jk}q^{np}a_i(m+p),
                 \end{equation}
                 \begin{equation}
                      [g_{ij}(m,n),a_k(p)a_l(s)]=0,
                 \end{equation}
                 \begin{equation}
                      [g_{ij}(m,n),a_k(p)a_l^*(s)]
                      =-\delta_{il}q^{-n(m+s)}a_k(p)a_j(m+s)+\rho \delta_{jl}q^{ns}a_k(p)a_i(m+s),
                 \end{equation}
                 \vspace{-0.5cm} \\
                 \parbox{1cm}{\begin{eqnarray}\end{eqnarray}}\hfill \parbox{13.66cm}
                    {\begin{eqnarray*}
                       [g_{ij}(m,n),a_k^*(p)a_l^*(s)]=-\delta_{ik}q^{-n(m+p)}a_j(m+p)a_l^*(s)+\rho
                       \delta_{jk}q^{np}a_i(m+p)a_l^*(s)&& \\-\delta_{il}q^{-n(m+s)}a_k^*(p)a_j(m+s)
                       +\rho \delta_{jl}q^{ns}a_k^*(p)a_i(m+s),&&  \hfill
                     \end{eqnarray*}}
                 \begin{equation}
                      [f_{ij}(m,n),a_k(p)]=\delta_{jk}q^{np}a_i(m+p),
                 \end{equation}
                 \begin{equation}
                      [f_{ij}(m,n),a_k^*(p)]=-\delta_{ik}q^{-n(m+p)}a_j^*(m+p),
                 \end{equation}
                 \begin{equation}
                      [f_{ij}(m,n),a_k(p)a_l^*(s)]
                      =\delta_{jk}q^{np}a_i(m+p)a_l^*(s)-\delta_{il}q^{-n(m+s)}a_k(p)a_j^*(m+p),
                 \end{equation}
                 \begin{equation}
                      [f_{ij}(m,n),a_k^*(p)a_l^*(s)]
                      =-\delta_{ik}q^{-n(m+p)}a_j^*(m+p)a_l^*(s)-\delta_{il}q^{-n(m+s)}a_k^*(p)a_j^*(m+s),
                 \end{equation}
                 \begin{equation}
                      [h_{ij}(m,n),a_k^*(p)]=0,
                 \end{equation}
                 \begin{equation}
                      [h_{ij}(m,n),a_k^*(p)a_l^*(s)]=0,
                 \end{equation}
         for m,n,p,s $\in \mathbb{Z}$ and $1 \leq i,j,k,l \leq N$.
            \end{lemma}

       \noindent \textbf{Proof.} First, we have
            \begin{eqnarray*}
                 &&[g_{ij}(m,n),a_k^*(p)]
                 \\
                 &=&\sum_{s \in \mathbb{Z}}q^{-ns}[:a_i(m-s)a_j(s):,a_k^*(p)]
                 \\
                 &=&\sum_{s \in \mathbb{Z}}q^{-ns}[a_i(m-s)a_j(s),a_k^*(p)]
                 \\
                 &=&\sum_{s \in \mathbb{Z}}q^{-ns} \Bigl(a_i(m-s)\bigl\{a_j(s),a_k^*(p)\bigr\}_{\rho}
                                        -\rho\bigl\{a_i(m-s),a_k^*(p)\bigr\}_{\rho}a_j(s)\Bigr)
                 \\
                 &=&-\delta_{ik}q^{-n(m+p)}a_j(m+p)+\rho
                 \delta_{jk}q^{np}a_i(m+p).
            \end{eqnarray*}
       Then
            \begin{eqnarray*}
                 &&[g_{ij}(m,n),a_k^*(p)a_l^*(s)]\\
                 &=&[g_{ij}(m,n),a_k^*(p)]a_l^*(s)+a_k^*(p)[g_{ij}(m,n),a_l^*(s)]   \\
                 &=&-\delta_{ik}q^{-n(m+p)}a_j(m+p)a_l^*(s)+\rho \delta_{jk}q^{np}a_i(m+p)a_l^*(s)\\
                 && \qquad -\delta_{il}q^{-n(m+s)}a_k^*(p)a_j(m+s)+\rho
                 \delta_{jl}q^{ns}a_k^*(p)a_i(m+s).
            \end{eqnarray*}
       So (2.18) and (2.21) hold true. The proof of the others is similar.             $\hfill \blacksquare$

       In what follows we shall mean $\frac{q^{mn}-1}{q^n-1}=m$ if $n \in \Lambda(q)$. This will make our
       formulas more concise.

       Next we list all Lie brackets that are needed.
            \begin{prop}
                 \begin{equation*}
                      [g_{ij}(m,n),g_{kl}(p,s)]=0
                 \end{equation*}
       for all $m,p,n,s \in \mathbb{Z}$ and $1 \leq i,j,k,l \leq N$.
            \end{prop}
            \begin{prop}
                 \begin{eqnarray*}
                      &&[g_{ij}(m,n),f_{kl}(p,s)]  \\
                      &=&-\delta_{il}q^{ms}g_{kj}(m+p,n+s)+\rho \delta_{jl}q^{(s-n)m}g_{ki}(m+p,s-n)
                 \end{eqnarray*}
       for all $m,p,n,s \in \mathbb{Z}$ and $1 \leq i,j,k,l \leq N$.
            \end{prop}
            \begin{prop}
               \begin{eqnarray*}
                  &&[g_{ij}(m,n),h_{kl}(p,s)]\\
                  &=&-\delta_{ik}q^{-n(m+p)}f_{jl}(m+p,s-n)+\rho \delta_{jk}q^{np}f_{il}(m+p,n+s)  \\
                  &&+\rho \delta_{il}q^{-(mn+np+ps)}f_{jk}(m+p,-(n+s))-\delta_{jl}q^{(n-s)p}f_{ik}(m+p,n-s) \\
                  &&-\rho \delta_{ik}\delta_{jl}\delta_{m+p,0}\frac{1}{2}(q^{s-n}+1)\frac{q^{m(s-n)}-1}{q^{s-n}-1}
                    +\delta_{jk}\delta_{il}\delta_{m+p,0}q^{np}\frac{1}{2}(q^{s+n}+1)\frac{q^{m(s+n)}-1}{q^{s+n}-1}
               \end{eqnarray*}
       for $m,p,n,s \in \mathbb{Z}$ and $1 \leq i,j,k,l \leq N$.
            \end{prop}
            \begin{prop}
                 \begin{eqnarray*}
                      [f_{ij}(m,n),f_{kl}(p,s)]
                       &=&\delta_{jk}q^{np}f_{il}(m+p,n+s)-\delta_{il}q^{sm}f_{kj}(m+p,n+s)\\
                       &&+\rho \delta_{jk} \delta_{il} q^{np} \delta_{m+p,0} \frac{1}{2}
                         (q^{s+n}+1)\frac{q^{m(s+n)}-1}{q^{s+n}-1}
                 \end{eqnarray*}
       for $m,p,n,s \in \mathbb{Z}$ and $1 \leq i,j,k,l \leq N$.
            \end{prop}
            \begin{prop}
               \begin{eqnarray*}
                  &[f_{ij}(m,n),h_{kl}(p,s)]=-\delta_{ik}q^{-n(m+p)}h_{jl}(m+p,s-n)-\delta_{il}q^{ms}h_{kj}(m+p,n+s)&
               \end{eqnarray*}
       for $m,p,n,s \in \mathbb{Z}$ and $1 \leq i,j,k,l \leq N$.
            \end{prop}
            \begin{prop}
                 \begin{equation*}
                      [h_{ij}(m,n),h_{kl}(p,s)]=0
                 \end{equation*}
       for $m,p,n,s \in \mathbb{Z}$ and $1 \leq i,j,k,l \leq N$.
            \end{prop}

            We shall only prove Proposition 2.3 which is the most
            complicated one. The proof of the others is either
            similar or easy.

       \noindent\textbf{Proof of Proposition 2.3 }It follows from (2.21) and (2.7) that
            \begin{eqnarray*}
                 &&[g_{ij}(m,n),q^{-st}:a_k^*(p-t)a_l^*(t):]  \\
                 &=&-\delta_{ik}q^{-st-n(m+p-t)}a_j(m+p-t)a_l^*(t)+\rho
                    \delta_{jk}q^{-st+n(p-t)}a_i(m+p-t)a_l^*(t) \\
                 &&-\delta_{il}q^{-st-n(m+t)}a_k^*(p-t)a_j(m+t)
                   +\rho \delta_{jl}q^{-st+nt}a_k^*(p-t)a_i(m+t) \\
                 &=&-\delta_{ik}q^{-st-n(m+p-t)}
                    \bigl(:a_j(m+p-t)a_l^*(t):+\rho \delta_{jl}\delta_{m+p,0}\theta(m+p-2t)\bigr)\\
                 &&+\rho \delta_{jk}q^{-st+n(p-t)}
                    \bigl(:a_i(m+p-t)a_l^*(t):+\rho \delta_{il}\delta_{m+p,0}\theta(m+p-2t)\bigr) \\
                 &&-\delta_{il}q^{-st-n(m+t)}
                   \bigl(-\rho:a_j(m+t)a_k^*(p-t):+\delta_{jk}\delta_{m+p,0}\theta(p-m-2t)\bigr)\\
                 &&+\rho \delta_{jl}q^{-st+nt}
                    \bigl(-\rho:a_i(m+t)a_k^*(p-t):+\delta_{ik}\delta_{m+p,0}\theta(p-m-2t)\bigr)\\
                 &=&-\delta_{ik}q^{-n(m+p)}q^{-(s-n)t}:a_j(m+p-t)a_l^*(t):           \\
                 && +\rho \delta_{jk}q^{np}q^{-(n+s)t}:a_i(m+p-t)a_l^*(t):       \\
                 &&+\rho\delta_{il}q^{-pn-ps-nm}q^{(s+n)(p-t)}:a_j(m+t)a_k^*(p-t):   \\
                 &&-\delta_{jl}q^{p(n-s)}q^{-(n-s)(p-t)}:a_i(m+t)a_k^*(p-t):     \\
                 &&-\rho\delta_{ik}\delta_{jl}\delta_{m+p,0}q^{-(s-n)t}(\theta(-2t)-\theta(-2m-2t)) \\
                 &&+\delta_{jk}\delta_{il}\delta_{m+p,0}q^{np}q^{-(n+s)t}(\theta(-2t)-\theta(-2m-2t)).
            \end{eqnarray*}
       Since
            \begin{eqnarray}
                 &&\sum_{t \in \mathbb{Z}}q^{-xt}\Bigl(\theta(-2t)-\theta(-2m-2t)\Bigr)   \nonumber  \\
                 &=&\left\{
                    \begin{array}
                       {l@{ \quad }l}
                          0, & \mbox{if } m=0,\vspace{0.1cm}  \\
                          \frac{1}{2}\bigl(1+q^{xm}\bigr)+\sum_{t=-(m-1)}^{-1}q^{-xt}, & \mbox{if }  m>0,  \\
                          -\frac{1}{2}\bigl(1+q^{xm}\bigr)-\sum_{t=1}^{-m-1}q^{-xt}, & \mbox{if } m<0
                    \end{array} \right.  \\
                 &=&\frac{q^{(m+1)x}-q^x+q^{mx}-1}{2(q^x-1)} \nonumber \\
                 &=&\frac{1}{2} (q^x+1)\frac{q^{mx}-1}{q^x-1},
                 \nonumber
            \end{eqnarray}
       we obtain Proposition 2.3.            $\hfill \blacksquare$

       Next we shall find the correspondence between
       $g_{ij}(m,n)$, $h_{ij}(m,n)$, $f_{ij}(m,n)$ and $\tilde{g}_{ij}(m,n)$,
       $\tilde{h}_{ij}(m,n)$, $\tilde{f}_{ij}(m,n)$. To this end, we
       have to modify our operators $g_{ij}(m,n)$, $h_{ij}(m,n)$,
       $f_{ij}(m,n)$.

       From Proposition 2.3, we see that , if $n+s \in \Lambda(q)$ and $n-s \in \Lambda(q)$,
            \begin{eqnarray*}
                  &&[g_{ij}(m,n),h_{kl}(p,s)]\\
                  &=&-\delta_{ik}q^{-n(m+p)}f_{jl}(m+p,s-n)+\rho \delta_{jk}q^{np}f_{il}(m+p,n+s)  \\
                  &&+\rho \delta_{il}q^{-(mn+np+ps)}f_{jk}(m+p,-(n+s))-\delta_{jl}q^{(n-s)p}f_{ik}(m+p,n-s) \\
                  &&-\rho \delta_{ik}\delta_{jl}\delta_{m+p,0}m
                    +\delta_{jk}\delta_{il}\delta_{m+p,0}q^{np}m.
            \end{eqnarray*}

       If $n+s \in \mathbb{Z} \setminus \Lambda(q)$ and $n-s \in
       \Lambda(q)$,
            \begin{eqnarray*}
                  &&[g_{ij}(m,n),h_{kl}(p,s)]\\
                  &=&-\delta_{ik}q^{-n(m+p)}f_{jl}(m+p,s-n)+\rho \delta_{jk}q^{np}f_{il}(m+p,n+s)  \\
                  &&+\rho \delta_{il}q^{-(mn+np+ps)}f_{jk}(m+p,-(n+s))-\delta_{jl}q^{(n-s)p}f_{ik}(m+p,n-s) \\
                  &&-\rho \delta_{ik}\delta_{jl}\delta_{m+p,0}m
                    +\delta_{jk}\delta_{il}\delta_{m+p,0}q^{np}\frac{1}{2}(q^{s+n}+1)\frac{q^{m(s+n)}-1}{q^{s+n}-1}  \\
                  &=&-\delta_{ik}q^{-n(m+p)}f_{jl}(m+p,s-n)-\delta_{jl}q^{(n-s)p}f_{ik}(m+p,n-s)                 \\
                  &&+\rho \delta_{jk}q^{np}\Bigl(f_{il}(m+p,n+s)
                       -\frac{\rho}{2}\delta_{il}\delta_{m+p,0}\frac{q^{n+s}+1}{q^{n+s}-1}\Bigr)   \\
                  &&+\rho\delta_{il}q^{-(mn+np+ps)}\Bigl(f_{jk}(m+p,-n-s)
                       -\frac{\rho}{2}\delta_{jk}\delta_{m+p,0}\frac{q^{-n-s}+1}{q^{-n-s}-1})\Bigr)  \\
                  &&-\rho \delta_{ik}\delta_{jl}\delta_{m+p,0}m.
            \end{eqnarray*}

       Similarly, if $n+s \in \Lambda(q)$ and $n-s \in \mathbb{Z} \setminus
       \Lambda(q)$,
            \begin{eqnarray*}
                  &&[g_{ij}(m,n),h_{kl}(p,s)]      \\
                  &=&\rho \delta_{jk}q^{np}f_{il}(m+p,n+s)+\rho \delta_{il}q^{-(mn+np+ps)}f_{jk}(m+p,-n-s)       \\
                  &&-\delta_{ik}q^{-n(m+p)}\Bigl(f_{jl}(m+p,s-n)
                       -\frac{\rho}{2}\delta_{jl}\delta_{m+p,0}\frac{q^{s-n}+1}{q^{s-n}-1}\Bigr)   \\
                  &&-\delta_{jl}q^{(n-s)p}\Bigl(f_{ik}(m+p,n-s)
                       -\frac{\rho}{2}\delta_{ik}\delta_{m+p,0}\frac{q^{n-s}+1}{q^{n-s}-1})\Bigr)  \\
                  &&+\delta_{jk}\delta_{il}\delta_{m+p,0}q^{np}m.
            \end{eqnarray*}

       By the above two relations, we have if $n+s,n-s \in \mathbb{Z} \setminus
       \Lambda(q)$,
            \begin{eqnarray*}
                  &&[g_{ij}(m,n),h_{kl}(p,s)]      \\
                  &=&\rho \delta_{jk}q^{np}f_{il}(m+p,n+s)+\rho \delta_{il}q^{-(mn+np+ps)}f_{jk}(m+p,-n-s)       \\
                  &&+\rho \delta_{jk}q^{np}\Bigl(f_{il}(m+p,n+s)
                       -\frac{\rho}{2}\delta_{il}\delta_{m+p,0}\frac{q^{n+s}+1}{q^{n+s}-1}\Bigr)   \\
                  &&+\rho\delta_{il}q^{-(mn+np+ps)}\Bigl(f_{jk}(m+p,-n-s)
                       -\frac{\rho}{2}\delta_{jk}\delta_{m+p,0}\frac{q^{-n-s}+1}{q^{-n-s}-1})\Bigr)  \\
                  &&-\delta_{ik}q^{-n(m+p)}\Bigl(f_{jl}(m+p,s-n)
                       -\frac{\rho}{2}\delta_{jl}\delta_{m+p,0}\frac{q^{s-n}+1}{q^{s-n}-1}\Bigr)   \\
                  &&-\delta_{jl}q^{(n-s)p}\Bigl(f_{ik}(m+p,n-s)
                       -\frac{\rho}{2}\delta_{ik}\delta_{m+p,0}\frac{q^{n-s}+1}{q^{n-s}-1})\Bigr).
            \end{eqnarray*}

       Using the same method, from Proposition 2.4 we have, if $n+s \in \Lambda(q)$,
            \begin{eqnarray*}
                 &&[f_{ij}(m,n),f_{kl}(p,s)]    \\
                 &=&\delta_{jk}q^{np}f_{il}(m+p,n+s)-\delta_{il}q^{sm}f_{kj}(m+p,n+s)
                       +\rho \delta_{jk} \delta_{il} q^{np} \delta_{m+p,0} m .
            \end{eqnarray*}

       If $n+s \in \mathbb{Z} \setminus \Lambda(q)$, then
            \begin{eqnarray*}
                 &&[f_{ij}(m,n),f_{kl}(p,s)]    \\
                 &=&\delta_{jk}q^{np}\Bigl(f_{il}(m+p,n+s)
                      -\frac{\rho}{2}\delta_{il}\delta_{m+p,0}\frac{q^{n+s}+1}{q^{n+s}-1}\Bigr)   \\
                 &&-\delta_{il}q^{sm}\Bigl(f_{kj}(m+p,n+s)
                      -\frac{\rho}{2}\delta_{jk}\delta_{m+p,0}\frac{q^{n+s}+1}{q^{n+s}-1})\Bigr).
            \end{eqnarray*}

       If we define \\
           \noindent \parbox{1cm}{\begin{eqnarray}\end{eqnarray}}\hfill \parbox{13.66cm}
               {\begin{eqnarray*}& F_{ij}(m,n)=\left\{
                   \begin{array}
                    {r@{ \quad }l} f_{ij}(m,n), & \mbox{  for } n \in \Lambda(q) \\
                    f_{ij}(m,n)-\frac{1}{2}\rho \delta_{ij}\delta_{m,0}\frac{q^n+1}{q^n-1},
                    &\mbox{  for } n \in \mathbb{Z}
                    \setminus \Lambda(q) \end{array} \right. & \\ \vspace{1cm}
                   &G_{ij}(m,n)=g_{ij}(m,n),\quad H_{ij}(m,n)=h_{ij}(m,n),&  \hfill
                \end{eqnarray*}}
       then we have

            \begin{theorem}
               $V(N,\rho)$ is a module for the Lie algebra $\widehat{\mathcal{G}_{\rho}}$ under the action given by
                \begin{eqnarray*}
                     \pi(\tilde{g}_{ij}(m,n))=G_{ij}(m,n), &&\qquad  \pi(\tilde{f}_{ij}(m,n))=F_{ij}(m,n),\\
                     \pi(\tilde{h}_{ij}(m,n))=H_{ij}(m,n),&&\qquad   \pi(c(n))=\frac{\rho}{2},\qquad \pi(c_y)=0.
                \end{eqnarray*}
            \end{theorem}

\subsection{ Type B}
       To consider BC$_N$-graded Lie algebras with grading
       subalgebra of type B$_N$,
        we require an extension of the algebra $\alpha (N,+)$. The generators
            \begin{equation}
                 \{e(m)|m \in \mathbb{Z}\}
            \end{equation}
       span an infinite-dimensional Cli\/f\/ford algebra with relations
            \begin{equation}
                 \{e(m),e(n)\}_{+}=e(m)e(n)+e(n)e(m)=\delta_{n+m,0}.
            \end{equation}
       Let $\alpha'(N)$ denote the algebra obtained by adjoining to $\alpha (N,+)$ the generators (2.30)
       with relations (2.31) and
            \begin{equation}
                 \{a_i(m),e(n)\}_{+}=0=\{a_i^*(m),e(n)\}_{+}
            \end{equation}

       We now define the normal ordering as in (2.6), i.e.  \\
            \parbox{1cm}{\begin{eqnarray}\end{eqnarray}}\hfill \parbox{13.66cm}
               {\begin{eqnarray*}& :e(m)e(n):=\left\{
                  \begin{array}
                        {r@{ }l}
                        & e(m)e(n) \hspace{3.18cm} \mbox{if }  n>m\vspace{0.1cm}
                        \\ &\frac{1}{2}\bigl(e(m)e(n)-e(n)e(m)\bigr) \quad \mbox{if }n=m
                        \\& -e(n)e(m) \hspace{2.88cm} \mbox{if }n<m
                  \end{array} \right., & \\
                            & :a_i(m)e(n):=a_i(m)e(n)=- e(n)a_i(m),& \\
                            & :a_i^*(m)e(n):=a_i^*(m)e(n)=- e(n)a_i^*(m),&  \hfill
                \end{eqnarray*}}
       for $n,m \in \mathbb{Z},1 \leq i,j \leq N.$ Then
            \begin{equation}
                 e(m)e(n)=:e(m)e(n):+\delta_{m+n,0} \theta(m-n).
            \end{equation}

       By (2.2), we have  \\
       \noindent \parbox{1cm}{\begin{eqnarray}\end{eqnarray}}\hfill \parbox{13.66cm}
                {\begin{eqnarray*}&& [a_i(m)a_i(n),e(p)]=[a_i(m)a_i^*(n),e(p)]=[a_i^*(m)a_i^*(n),e(p)]=0,  \\
                                  && [a_i(m)e(n),a_k(p)]=0,  \\
                                  && [a_i(m)e(n),a_k^*(p)]=-\delta_{ik} \delta_{m+p,0} e(n),  \\
                                  && [a_i(m)e(n),e(p)]=\delta_{n+p,0} a_i(m),  \\
                                  && [a_i^*(m)e(n),a_k^*(p)]=0, \\
                                  && [a_i^*(m)e(n),e(p)]=\delta_{n+p,0} a_i^*(m).\\
                                  && [e(m)e(n),e(p)]=\delta_{n+p,0}e(m)-\delta_{m+p,0}e(n).  \hfill
                \end{eqnarray*}}
       for $m,n,p \in \mathbb{Z},1 \leq i,j,k \leq N.$

       Let $V_0$ be a simple Cli\/f\/ford module for the Cli\/f\/ford algebra generated by (2.30) with relations
       (2.31) and containing``vacuum vector" $v'_0$, which is killed by annihilation operators.
       (Here we call $e(m)$
       annihilation operator if $m>0$, or a creation operator if $m<0$. e(0) acts as
       scalar.)Because of (2.32), we see that the $\alpha'(N)$-module
            \begin{equation}
                 V'(N)=V(N,+) \otimes V_0=\alpha'(N)v'_0
            \end{equation}
       is simple.

       Now we construct a class of fermions on $V'(N)$. For any $m,n \in \mathbb{Z},1 \leq i,j \leq N,$ set
            \begin{equation}
                 f_{ij}(m,n)=\sum_{s \in \mathbb{Z}} q^{-ns}:a_i(m-s)a_j^*(s):,
            \end{equation}
            \begin{equation}
                 g_{ij}(m,n)=\sum_{s \in \mathbb{Z}} q^{-ns}:a_i(m-s)a_j(s):,
            \end{equation}
            \begin{equation}
                 h_{ij}(m,n)=\sum_{s \in \mathbb{Z}} q^{-ns}:a_i^*(m-s)a_j^*(s):,
            \end{equation}
            \begin{equation}
                 e_i(m,n)=\sum_{s \in \mathbb{Z}} q^{-ns}:a_i(m-s)e(s):,
            \end{equation}
            \begin{equation}
                 e_i^*(m,n)=\sum_{s \in \mathbb{Z}} q^{-ns}:a_i^*(m-s)e(s):,
            \end{equation}
            \begin{equation}
                 e_0(m,n)=\sum_{s \in \mathbb{Z}}
                 q^{-ns}:e(m-s)e(s):.
            \end{equation}
            \begin{rem}
                 In this section, $g_{ij}(m,n),f_{ij}(m,n),h_{ij}(m,n)$ are
                 the same as ones in the type D case (2.13)-(2.15)
                  by taking $\rho=1$.  So we needn't to check the Lie brackets among them.
            \end{rem}
            \begin{lemma}
               We have
                 \begin{equation}
                      [g_{ij}(m,n),a_k(p)e(s)]=[g_{ij}(m,n),e(p)e(s)]=0,
                 \end{equation}
                 \begin{equation}
                      [g_{ij}(m,n),a_k^*(p)e(s)]=-\delta_{ik}q^{-n(m+p)}a_j(m+p)e(s)+ \delta_{jk}q^{np}a_i(m+p)e(s),
                 \end{equation}
                 \begin{equation}
                      [f_{ij}(m,n),a_k(p)e(s)]=\delta_{jk}q^{np}a_i(m+p)e(s),
                 \end{equation}
                 \begin{equation}
                      [f_{ij}(m,n),a_k^*(p)e(s)]=-\delta_{ik}q^{-n(m+p)}a_j^*(m+p)e(s),
                 \end{equation}
                 \begin{equation}
                      [f_{ij}(m,n),e(p)e(s)]=0,
                 \end{equation}
                 \begin{equation}
                      [h_{ij}(m,n),a_k(p)e(s)]=\delta_{jk}q^{np}a_i^*(m+p)e(s)-
                      \delta_{ik}q^{-n(m+p)}a_j^*(m+p)e(s),
                 \end{equation}
                 \begin{equation}
                      [h_{ij}(m,n),a_k^*(p)e(s)]=[h_{ij}(m,n),e(p)e(s)]=0,
                 \end{equation}
                 \begin{equation}
                      [e_i(m,n),a_k(p)]=0,
                 \end{equation}
                 \begin{equation}
                      [e_i(m,n),a_k^*(p)]=-\delta_{ik}q^{-n(m+p)}e(m+p),
                 \end{equation}
                 \begin{equation}
                      [e_i(m,n),e(p)]=q^{np}a_i(m+p),
                 \end{equation}
                 \begin{equation}
                      [e_i(m,n),a_k(p)e(s)]=q^{ns}a_k(p)a_i(m+s),
                 \end{equation}
                 \begin{equation}
                      [e_i(m,n),a_k^*(p)e(s)]=-\delta_{ik}q^{-n(m+p)}e(m+p)e(s)+q^{ns}a_k^*(p)a_i(m+s),
                 \end{equation}
                 \begin{equation}
                      [e_i(m,n),e(p)e(s)]=q^{np}a_i(m+p)e(s)+q^{ns}a_i(m+s)e(p),
                 \end{equation}
                 \begin{equation}
                      [e_i^*(m,n),a_k^*(p)]=0,
                 \end{equation}
                 \begin{equation}
                      [e_i^*(m,n),e(p)]=q^{np}a_i^*(m+p),
                 \end{equation}
                 \begin{equation}
                      [e_i^*(m,n),a_k^*(p)e(s)]=q^{ns}a_k^*(p)a_i^*(m+s),
                 \end{equation}
                 \begin{equation}
                      [e_i^*(m,n),e(p)e(s)]=q^{np}a_i^*(m+p)e(s)+q^{ns}a_i^*(m+s)e(p),
                 \end{equation}
                 \begin{equation}
                      [e_0(m,n),e(p)]= (q^{np}-q^{-n(m+p)})e(m+p),
                 \end{equation}
                 \begin{equation}
                      [e_0(m,n),e(p)e(s)]=(q^{np}-q^{-n(m+p)})e(m+p)e(s)+(q^{ns}-q^{-n(m+s)})e(p)e(m+s),
                 \end{equation}
               for m,n,p,s $\in \mathbb{Z}$ and $1 \leq i,j,k \leq N$.
            \end{lemma}
       As in Section 2.1, we have Propositions 2.1-2.6  plus the following
       propositions.

            \begin{prop}
                 \begin{eqnarray*}
                      &[g_{ij}(m,n),e_k(p,s)]=[g_{ij}(m,n),e_0(p,s)]=0,&\\
                      &[g_{ij}(m,n),e_k^*(p,s)]=-\delta_{ik}q^{-n(m+p)}e_j(m+p,s-n)+ \delta_{jk}q^{np}e_i(m+p,n+s)&
                 \end{eqnarray*}
               for all $m,p,n,s \in \mathbb{Z}$ and $1 \leq i,j,k \leq N$.
            \end{prop}
            \begin{prop}
                 \begin{eqnarray*}
                      &[f_{ij}(m,n),e_k(p,s)]=\delta_{jk}q^{np}e_i(m+p,n+s),&\\
                      &[f_{ij}(m,n),e_k^*(p,s)]=-\delta_{ik}q^{-n(m+p)}e_j^*(m+p,s-n),&  \\
                      &[f_{ij}(m,n),e_0(p,s)]=0&
                 \end{eqnarray*}
               for all $m,p,n,s \in \mathbb{Z}$ and $1 \leq i,j,k \leq N$.
            \end{prop}
            \begin{prop}
                 \begin{eqnarray*}
                      &[h_{ij}(m,n),e_k(p,s)]=\delta_{jk}q^{np}e_i^*(m+p,n+s)- \delta_{ik}q^{-n(m+p)}e_j^*(m+p,s-n),&\\
                      & [h_{ij}(m,n),e_k^*(p,s)]=[h_{ij}(m,n),e_0(p,s)]=0&
                 \end{eqnarray*}
               for all $m,p,n,s \in \mathbb{Z}$ and $1 \leq i,j,k \leq N$.
            \end{prop}
            \begin{prop}
                 \begin{eqnarray*}
                      &&[e_i(m,n),e_k(p,s)]=q^{m(s-n)}g_{ki}(m+p,s-n),\\
                      &&[e_i(m,n),e_k^*(p,s)]=-\delta_{ik}q^{-n(m+p)}e_0(m+p,s-n)-q^{p(n-s)}f_{ik}(m+p,n-s) \\
                      &&\hspace{4cm}-\delta_{ik}\delta_{m+p,0}\frac{1}{2}(q^{s-n}+1)\frac{q^{m(s-n)}-1}{q^{s-n}-1}, \\
                      && [e_i(m,n),e_0(p,s)]=q^{np}e_i(m+p,n+s)-q^{p(n-s)}e_i(m+p,n-s)
                 \end{eqnarray*}
               for all $m,p,n,s \in \mathbb{Z}$ and $1 \leq i,k \leq N$.
            \end{prop}
            \begin{prop}
                 \begin{eqnarray*}
                      && [e_i^*(m,n),e_k^*(p,s)]=q^{m(s-n)}h_{ki}(m+p,s-n),  \\
                      && [e_i^*(m,n),e_0(p,s)]=q^{np}e_i^*(m+p,n+s)-q^{p(n-s)}e_i^*(m+p,n-s)
                 \end{eqnarray*}
               for all $m,p,n,s \in \mathbb{Z}$ and $1 \leq i,k \leq N$.
            \end{prop}
            \begin{prop}
                 \begin{eqnarray*}
                      &&[e_0(m,n),e_0(p,s)]    \\
                      &=&(q^{np}-q^{sm})e_0(m+p,n+s)
                          +\delta_{m+p,0}q^{np}\frac{1}{2}(q^{n+s}+1)\frac{q^{m(n+s)}-1}{q^{n+s}-1} \\
                      &&+(q^{m(s-n)}-q^{-n(m+p)})e_0(m+p,s-n)-
                         \delta_{m+p,0}\frac{1}{2}(q^{s-n}+1)\frac{q^{m(s-n)}-1}{q^{s-n}-1}
                 \end{eqnarray*}
               for all $m,p,n,s \in \mathbb{Z}$.
            \end{prop}

                 We only
                 give proofs for Proposition 2.10 and Proposition 2.12.
                  The proof for the others is either similar or easy.

       \noindent\textbf{Proof of Proposition 2.10 and Proposition 2.12.} \\

       First, it follows from (2.53)-(2.55), (2.34) and (2.7) that
            \begin{eqnarray*}
                 &&[e_i(m,n),q^{-st}:a_k(p-t)e(t):]=q^{nt-st}a_k(p-t)a_i(m+t)        \\
                 &=&q^{m(s-n)}q^{-(s-n)(m+t)}a_k(p-t)a_i(m+t)   \\
                 &=&q^{m(s-n)}q^{-(s-n)(m+t)}:a_k(p-t)a_i(m+t):,
            \end{eqnarray*}
            \begin{eqnarray*}
                 &&[e_i(m,n),q^{-st}:a_k^*(p-t)e(t):]                                \\
                 &=&q^{-st}\Bigl(-\delta_{ik}q^{-n(m+p-t)}e(m+p-t)e(t)-q^{nt}a_k^*(p-t)a_i(m+t)\Bigr)              \\
                 &=&-\delta_{ik}q^{-n(m+p)}q^{-(s-n)t}e(m+p-t)e(t)+q^{-(s-n)t}a_k^*(p-t)a_i(m+t)                   \\
                 &=&-\delta_{ik}q^{-n(m+p)}q^{-(s-n)t}\Bigl(:e(m+p-t)e(t):+\delta_{m+p,0}\theta(m+p-2t)\Bigr)      \\
                 &&-q^{-(s-n)t}\bigl(:a_i(m+t)a_k^*(p-t):-\delta_{ik}\delta_{m+p,0}\theta(p-m-2t)\bigr)            \\
                 &=&-\delta_{ik}q^{-n(m+p)}q^{-(s-n)t}:e(m+p-t)e(t):   \\
                   &&-q^{p(n-s)}q^{-(n-s)(p-t)}:a_i(m+t)a_k^*(p-t): \\
                 &&-\delta_{ik}\delta_{m+p,0}q^{-(s-n)t}\bigl(\theta(-2t)-\theta(-2m-2t)
                 \bigr),
            \end{eqnarray*}
            \begin{eqnarray*}
                 &&[e_i(m,n),q^{-st}:e(p-t)e(t):]      \\
                 &=&q^{-st}\bigl(q^{n(p-t)}a_i(m+p-t)e(t)+q^{nt}a_i(m+t)e(p-t)\bigr)      \\
                 &=&q^{np}q^{-(n+s)t}:a_i(m+p-t)e(t):+q^{p(n-s)}q^{-(n-s)(p-t)}:a_i(m+t)e(p-t):.     \\
            \end{eqnarray*}
       Then by (2.28), we see that Proposition 2.10 holds true. \\

       Secondly, it follows from (2.61), (2.34) and (2.28) that
            \begin{eqnarray*}
                 &&[e_0(m,n),q^{-st}:e(p-t)e(t):]    \\
                 &=&q^{-st}\Bigl((q^{n(p-t)}-q^{-n(m+p-t)})e(m+p-t)e(t)+(q^{nt}-q^{-n(m+t)}) e(p-t)e(m+t)\Bigr) \\
                 &=&q^{-st}(q^{n(p-t)}-q^{-n(m+p-t)})\bigl(:e(m+p-t)e(t):+\delta_{m+p,0}\theta(m+p-2t)\bigr)\\
                 &&+q^{-st}(q^{nt}-q^{-n(m+t)})\bigl(:e(p-t)e(m+t):+\delta_{m+p,0}\theta(p-m-2t)\bigr) \\
                 &=&q^{np}q^{-(s+n)t}:e(m+p-t)e(t):-q^{-n(m+p)}q^{-(s-n)t}:e(m+p-t)e(t):\\
                  &&+q^{m(s-n)}q^{-(s-t)(m+t)}:e(p-t)e(m+t):-q^{sm}q^{-(n+s)(m+t)}:e(p-t)e(m+t): \\
                  &&+\delta_{m+p,0}q^{np}q^{-(n+s)t}\bigl(\theta(-2t)-\theta(-2m-2t)\bigr)  \\
                  &&-\delta_{m+p,0}q^{-(s-n)t}\bigl(\theta(-2t)-\theta(-2m-2t)\bigr)
            \end{eqnarray*}
       and Proposition 2.12 holds true.       $\hfill \blacksquare$

       As in Section 2.1 of type D case, we need to modify the
       definition of our operators.

       For Proposition 2.10, if $n-s \in \Lambda(q)$,
            \begin{eqnarray*}
                 &&[e_i(m,n),e_k^*(p,s)]    \\
                 &=&-\delta_{ik}q^{-n(m+p)}e_0(m+p,s-n)-q^{p(n-s)}f_{ik}(m+p,n-s)
                      -\delta_{ik}\delta_{m+p,0}m;
            \end{eqnarray*}

       if $n-s \in \mathbb{Z} \setminus \Lambda(q)$,
            \begin{eqnarray*}
                 &&[e_i(m,n),e_k^*(p,s)]    \\
                 &=&-\delta_{ik}q^{-n(m+p)}\Bigl(e_0(m+p,s-n)
                      -\frac{1}{2}\delta_{m+p,0}\frac{q^{s-n}+1}{q^{s-n}-1}\Bigr) \\
                 &&-q^{p(n-s)}\Bigl(f_{ik}(m+p,n-s)
                      -\frac{1}{2}\delta_{jk}\delta_{m+p,0}\frac{q^{n-s}+1}{q^{n-s}-1})\Bigr).
            \end{eqnarray*}

       For Proposition 2.12, if $n+s \in \Lambda(q)$ and $n-s \in \Lambda(q)$,
            \begin{eqnarray*}
                 &&[e_0(m,n),e_0(p,s)]    \\
                 &=&(q^{np}-q^{sm})e_0(m+p,n+s)
                          +(q^{m(s-n)}-q^{-n(m+p)})e_0(m+p,s-n) \\
                 &&+\delta_{m+p,0}q^{np}m-\delta_{m+p,0}m;
            \end{eqnarray*}

       if $n+s \in \mathbb{Z} \setminus \Lambda(q)$ and $n-s \in
       \Lambda(q)$,
            \begin{eqnarray*}
                  &&[e_0(m,n),e_0(p,s)]    \\
                      &=&(q^{np}-q^{sm})\Bigl(e_0(m+p,n+s)
                          -\frac{1}{2}\delta_{m+p,0}\frac{q^{n+s}+1}{q^{n+s}-1}\Bigr) \\
                      &&+(q^{m(s-n)}-q^{-n(m+p)})e_0(m+p,s-n)-\delta_{m+p,0}m;
            \end{eqnarray*}

       if $n+s \in \Lambda(q)$ and $n-s \in \mathbb{Z} \setminus
       \Lambda(q)$,
            \begin{eqnarray*}
                  &&[e_0(m,n),e_0(p,s)]    \\
                  &=&(q^{np}-q^{sm})e_0(m+p,n+s)+\delta_{m+p,0}q^{np}m                \\
                  &&+(q^{m(s-n)}-q^{-n(m+p)})\Bigl(e_0(m+p,s-n)
                         -\frac{1}{2}\delta_{m+p,0}\frac{q^{s-n}+1}{q^{s-n}-1}\Bigr);
            \end{eqnarray*}

       if $n+s,n-s \in \mathbb{Z} \setminus \Lambda(q)$,
            \begin{eqnarray*}
                  &&[e_0(m,n),e_0(p,s)]    \\
                  &=&(q^{np}-q^{sm})\Bigl(e_0(m+p,n+s)
                         -\frac{1}{2}\delta_{m+p,0}\frac{q^{n+s}+1}{q^{n+s}-1}\Bigr)           \\
                  &&+(q^{m(s-n)}-q^{-n(m+p)})\Bigl(e_0(m+p,s-n)
                         -\frac{1}{2}\delta_{m+p,0}\frac{q^{s-n}+1}{q^{s-n}-1}\Bigr).
            \end{eqnarray*}

       Now we define \\
            \noindent \parbox{1cm}{\begin{eqnarray}\end{eqnarray}}\hfill \parbox{13.66cm}
               {\begin{eqnarray*}& F_{ij}(m,n)=\left\{
                   \begin{array}
                       {r@{ \quad }l} f_{ij}(m,n), & \mbox{  for } n \in \Lambda(q),  \\
                       f_{ij}(m,n)-\frac{1}{2}\delta_{ij}\delta_{m,0}\frac{q^n+1}{q^n-1},
                          &\mbox{  for } n \in \mathbb{Z}   \setminus \Lambda(q),
                   \end{array} \right. & \\ \vspace{1cm}
                       &G_{ij}(m,n)=g_{ij}(m,n),\quad H_{ij}(m,n)=h_{ij}(m,n),& \\
                       &E_i(m,n)=e_i(m,n),\quad E_i^*(m,n)=e_i^*(m,n),& \\
                       & E_0(m,n)=\left\{
                   \begin{array}
                        {r@{ \quad }l} e_0(m,n), & \mbox{  for } n \in \Lambda(q), \\
                        e_0(m,n)-\frac{1}{2}\delta_{m,0}\frac{q^n+1}{q^n-1}, &\mbox{  for } n \in \mathbb{Z}
                        \setminus \Lambda(q).
                   \end{array} \right. &   \hfill
                \end{eqnarray*}}
       Then we have

            \begin{theorem}
                $V'(N)$ is a module for the Lie algebra $\widehat{\mathcal{G'}}$ under the action given by
                \begin{eqnarray*}
                     \pi(\tilde{g}_{ij}(m,n))=G_{ij}(m,n), &&\qquad  \pi(\tilde{f}_{ij}(m,n))=F_{ij}(m,n),\\
                     \pi(\tilde{h}_{ij}(m,n))=H_{ij}(m,n), &&\qquad  \pi(\tilde{e}_i(m,n))=E_i(m,n),      \\
                     \pi(\tilde{e}_i^*(m,n))=E_i^*(m,n),   &&\qquad  \pi(\tilde{e}_0(m,n))=E_0(m,n),      \\
                     \pi(c(n))=\frac{1}{2},                &&\qquad  \pi(c_y)=0.
                \end{eqnarray*}
            \end{theorem}

       \footnotesize{    \begin{align*}
%         &\mbox{Hongjia Chen}                                    &&\mbox{Yun Gao}          \\
         &\mbox{Department of Mathematics}                       &&\mbox{Department of Mathematics and Statistics} \\
         &\mbox{University of Science and Technology of China}   &&\mbox{York University} \\
         &\mbox{Hefei, Anhui}                                    &&\mbox{Toronto, Ontario} \\
         &\mbox{P. R. China  230026}                             &&\mbox{Canada  M3J 1P3} \\
         &\mbox{Email:hjchen@mail.ustc.edu.cn}                   &&\mbox{Email:ygao@yorku.ca} \\
       \end{align*}}

\end{document}